\documentclass[12pt]{article}   

\usepackage{graphicx,amssymb,bbm,amsmath,mathrsfs,amsthm,amsfonts}
\usepackage[margin=1in]{geometry}        
\usepackage{amsfonts}             
\usepackage{amsthm}

\newtheorem{theorem}{Theorem}[section]
\newtheorem{lemma}[theorem]{Lemma}

\newtheorem{proposition}[theorem]{Proposition}

\newtheorem{remark}[theorem]{Remark}
\newtheorem{example}[theorem]{Example}

\begin{document}
	
\nocite{*}	

\title{The structure of entrance laws for time-inhomogeneous Ornstein-Uhlenbeck Processes with L\'evy Noise in Hilbert spaces}

\author{Narges Rezvani Majid \\
              Bielefeld University \\
              rezvani.narges@gmail.com \\
           \and
           Michael R\"ockner \\
           	  Bielefeld University \&\\
           	   AMSS, CAS, Beijing\\         
           	  roeckner@math.uni-bielefeld.de   
}

\maketitle

\begin{abstract}
This paper is about the structure of all entrance laws (in the sense of Dynkin) for time-inhomogeneous Ornstein-Uhlenbeck processes with L\'evy 
noise in Hilbert state spaces. We identify the extremal entrance laws with finite weak first moments through an explicit formula for their Fourier transforms, generalising corresponding results by Dynkin for Wiener noise and nuclear state spaces. We then prove that an arbitrary entrance law with finite weak first moments can be uniquely represented as an integral over extremals. It is proved that this can be derived from Dynkin's seminal work "Sufficient statistics and extreme points" in Ann. Probab. 1978, which contains a purely measure theoretic generalization of the classical analytic Krein-Milman and Choquet Theorems. As an application, we obtain an easy uniqueness proof for $T$-periodic entrance laws in the general periodic case. A number of further applications to concrete cases are presented. 
\end{abstract}	
\\
\\
Dedicated to the memory of E. B. Dynkin.

\section{Introduction}
Let $H$ be a separable Hilbert space with Borel $\sigma$-algebra $\mathscr{B}(H)$. Consider a Markovian family of transition probabilities $\pi=(\pi_{s,t})_{s\leq t}$, i.e.,
\begin{description}
	\item[{\emph{(i)}}] $\pi_{s,t}(x,\cdot)$ is a probability measure on $(H,\mathscr{B}(H))$ for each $s\leq t$, $x\in H$.
	\item[{\emph{(ii)}}] $\pi_{s,t}(\cdot,B)$ belongs to $\mathcal{B}_b(H)$ (:= the set of all real-valued bounded measurable functions on $H$), for each $s\leq t$, $B\in\mathscr{B}(H)$.
	\item[{\emph{(iii)}}] $\pi_{s,t}(x,B)=\int_H \pi_{s,r}(x,dy)\pi_{r,t}(y,B)$, for each $s\leq r\leq t$, $B\in\mathscr{B}(H)$.
	\item[{\emph{(iv)}}] $\pi_{s,s}(x,B)=\mathbbm{1}_B(x)$, for each $x\in H$, $B\in \mathscr{B}(H)$.
\end{description}

Typical examples of such families are the transition probabilities of solutions to stochastic differential equations, whose drift and diffusion coefficients are 
time-dependent, but not random.

In this paper, we study {\emph{entrance laws}} or {\emph{evolution systems of measures}} corresponding to such transition probabilities  
introduced by E. B. Dynkin in \cite{MR0518321}. 
They are defined as families of probability measures $(\nu_t)_{t\in\mathbb{R}}$ on $H$ such that for all $s\leq t$,
\begin{align*}
\int_H\pi_{s,t}(x,B)\nu_s(dx)=\nu_{t}(B),\quad s\leq t,\ s,t\in\mathbb{R},\quad B\in\mathscr{B}(H)
\end{align*}
or in short
\begin{align*}
\nu_s \pi_{s,t}=\nu_t,\quad s\leq t.
\end{align*}
For example, if $\pi$ is time homogeneous, i.e., $\pi_{s,t}=\pi_{0,t-s}$, $s\leq t$, and has an invariant measure $\nu$, then $\nu_t:=\nu$, $t\in\mathbb{R}$, is a 
particular case of an entrance law.\\
We denote the set of all probability $\pi$-entrance laws $(\nu_t)_{t\in\mathbb{R}}$ by $\mathcal{K}(\pi)$. Obviously, $\mathcal{K}(\pi)$ forms a convex set and generically
it consists
of more than one element.

In his seminal work \cite{MR0518321} E. B. Dynkin proved a purely measure theoretic analogue of the corresponding well-known analytic results by Choquet or Krein and Milman, which
states that $\mathcal{K}(\pi)$ is a simplex, i.e., each element in $\mathcal{K}(\pi)$ has a unique integral representation
in terms of its extreme points.
Therefore, to fully understand the structure of $\mathcal{K}(\pi)$, it suffices to characterize the set of all extreme points denoted by $\mathcal{K}_e(\pi)$.\\

We concentrate on an important class of Markovian families of transition \break probabilities 
$(\pi_{s,t})_{s\leq t}$, associated with time-inhomogeneous Ornstein-Uhlenbeck processes with L\'evy noise. 
In the time homogeneous case such Ornstein-Uhlenbeck processes ``with jumps'' and their corresponding transition semigroups, called generalized Mehler semigroups, have been studied intensively, 
see \cite{MR2243880,MR2048508,MR1930955,MR2032114,MR2760602,MR2105913,Lunardi2019SchauderTF,MR3466574,MR2886463,MR1996872,MR1841407}.
In this paper, however, we look at the more general time-inhomogeneous case: \\
Let $(A(t),\mathcal{D}(A(t)))_{t\in\mathbb{R}}$ be a family of linear operators on $H$ with dense domains. Suppose that the non-autonomous Cauchy problem 
\begin{align*}
dX(t)=A(t)X(t)dt,\quad X(s)=x\in\mathcal{D}(A(s)),\quad s\leq t,
\end{align*}
is well-posed in the mild sense and has a unique solution given by a strong evolution family of linear operators $(U_{s,t})_{s\leq t}$ on $H$, where here and below $s,t$ run through all of $\mathbb{R}$. Recall  that $U=(U_{s,t})_{s\leq t}$ is a 
strong evolution family of bounded linear operators on $H$, if each  $U_{s,t}\in\mathcal{L}(H)$, $U_{t,t} = I$ for all $t\in\mathbb{R}$, $U_{r,t}U_{s,r} = U_{s,t}$ for 
all $s\leq r\leq t$ and $U$ is strongly continuous on $\{(s,t)\in\mathbb{R}^2\mid s\leq t\}$. Here $\mathcal{L}(H)$ denotes the set of all bounded linear operators on $H$.\\
We consider the following type of stochastic differential equations on $H$:
\begin{equation}\label{2'}
\begin{array}{ll}
dX(t)=A(t)X(t)dt+\sigma(t)dL(t),\quad s\leq t,\\
\ X(s)=x,
\end{array}
\end{equation}
where $\sigma:\mathbb{R}\rightarrow\mathcal{L}(H)$ is strongly measurable and $L$ is an 
$H$-valued L\'evy process.\\
Let $X(s,t,x)$, $s\leq t$, be the mild solution of equation $(\ref{2'})$, i.e.
\begin{align}\label{R1}
X(s,t,x)=U_{s,t}x+\int_s^t U_{r,t}\sigma(r)dL(r),\quad s\leq t,\ x\in H.
\end{align}
This mild solution is called time-inhomogeneous Ornstein-Uhlenbeck process with L\'evy noise.
Then, the associated family of transition probabilities $\pi=(\pi_{s,t})_{s\leq t}$ is called a time-inhomogeneous (generalized) Mehler semigroup, which is defined by:
\begin{align}\label{R2}
\pi_{s,t}(x,dy)=\mathbb{P}\circ X(s,t,x)^{-1}\ (dy)=\mu_{s,t}(dy-U_{s,t}x),\quad s\leq t, 
\end{align}
where $\mu_{s,t}$ is the distribution of the stochastic convolution $\int_s^tU_{r,t}\sigma(r)dL(r)$.\\
Generalized Mehler semigroups were initially
defined by Bogachev, R\"ockner, and Schmuland \cite{MR1392452} in the case of Wiener noise. This was extended to the non-Gaussian case in \cite{MR1745332}. The time inhomogeneous non-Gaussian case was studied in \cite{MR2861314} and further generalized in \cite{MR3466574}.

The question whether $(\ref{2'})$ has a solution in the sense of $(\ref{R1})$ reduces to the question whether the stochastic integral in $(\ref{R1})$ makes sense. In this respect we refer to \cite{MR2861314}, because in this paper we shall solely concentrate on the Markovian transition probabilities in $(\ref{R2})$, so only need the existence of the measures $\mu_{s,t}$, $s\leq t$.

Let $\mathcal{K}^1(\pi)$ be the set of all elements of $\mathcal{K}(\pi)$ with finite first weak moments and let $\mathcal{K}(U)$ denote the set of all 
$\kappa=(\kappa_t)_{t\in\mathbb{R}}\subset{H}$ with $U_{s,t}\kappa_s=\kappa_{t}$ for all $s\leq t$. Then
the main result of this paper (Theorem \ref{27!}), states, that under wide conditions, there exists a one-to-one 
correspondence between $\mathcal{K}(U)$ and the set $\mathcal{K}^1_e(\pi)$ of all extremal points of $\mathcal{K}^1(\pi)$. Furthermore, 
we show that the extremal $\pi$-entrance laws have explicit characteristic functions of the form $(\ref{09!})$ below. Moreover, we show that $\mathcal{K}^1(\pi)$ is a simplex (see Theorem \ref{27!}).

In the particular case of time-inhomogeneous Ornstein-Uhlenbeck processes with Wiener noise, a similar result was obtained by E. B. Dynkin in \cite{MR1013934} (see  Theorem 5.1 in there), however, with a family of nuclear spaces replacing our Hilbert space $H$ and assuming that such nuclear spaces exist satisfying all properties used for the proof. We generalize this result to time-inhomogeneous
Ornstein-Uhlenbeck processes with
L\'evy noise and implement this in a Hilbert space setting giving explicit (checkable) wide conditions under which our result holds.

This paper is organized as follows. In Section 2 we construct the time-\break inhomogeneous Mehler semigroups by using their characteristic
functions.
Section 3 is the main part of this paper, where the explicit formula for the characteristic functions of the extremal $\pi$-entrance laws is derived. This result is stated in 
Theorem \ref{27!}. In Section \ref{S4}, we will show how Theorem \ref{27!} can be applied to prove
uniqueness of ($T$-periodic) $\pi$-entrance laws (see Theorem \ref{3'}). Section 5 is devoted to examples. .

\section{Definitions, hypotheses and construction}
Let us fix a real separable Hilbert space $H$ with inner product $\langle\cdot,\cdot\rangle$ and corresponding norm $\lVert\cdot\rVert$. For a probability measure $\mu$ on $(H,\mathscr{B}({H}))$, 
we recall that its characteristic
function is defined by 
\begin{equation*}
\widehat{\mu}(a)=\int_{H}e^{i\langle a,x\rangle}\mu(dx),\quad a\in{H}.
\end{equation*}
We recall that by a monotone class argument, every probability measure $\mu$ is uniquely determined by its characteristic function $\widehat{\mu}$.\\

We also recall that a function $\varphi:{H}\rightarrow\mathbb{C}$ is called positive definite if for all $n\in\mathbb{N},\ a_1,...,a_n\in{H}$ and $c_1,...,c_n\in\mathbb{C}$,
\begin{equation*}
\sum_{i,j=1}^{n}\varphi(a_i-a_j)c_i\overline{c_j}\geq0
\end{equation*}
and $\varphi$ is called negative definite if $\varphi(0)\geq 0$, $\varphi(-a)=\overline{\varphi(a)}$ for all $a\in H$ and 
for all $n\in\mathbb{N},\ a_1,...,a_n\in{H}$ and $c_1,...,c_n\in\mathbb{C}$
with $\sum_{i=1}^{n}c_i=0$, we have 
\begin{equation*}
\sum_{i,j=1}^{n}\varphi(a_i-a_j)c_i\overline{c_j}\leq0.
\end{equation*}

The \emph{Sazonov topology} is the topology on $H$ generated by the set of seminorms $a\mapsto\lVert Sa\rVert$, $a\in H$, where $S$ ranges over 
the family of all Hilbert-Schmidt operators on $H$.

By the {\it{Minlos-Sazonov theorem}} (see e.g. Theorem 2.4, Chapter VI in \cite{MR0226684} or Theorem VI.1.1 in \cite{MR1435288}), a complex-valued 
function $\varphi$ on ${H}$, is the characteristic function
of a probability measure on $(H,\mathscr{B}(H))$ if and only if

${\emph{(i)}}\ \varphi(0)=1$,

${\emph{(ii)}}\ \varphi$ is positive definite on $H$,

${\emph{(iii)}}\ \varphi$ is Sazonov continuous on $H$.\\

Let $\mathcal{L}_1^+(H)$ denote the set of all non-negative symmetric trace class operators on $H$, which is a Banach space with norm $\lVert \cdot\rVert_{\mathcal{L}_1^+}$.
By the {\emph{L\'evy-Khinchin formula}}, a function $\varphi:H\rightarrow\mathbb{C}$ is the characteristic function of an infinitely divisible probability measure
$\mu$ (see Definition 4.1 in Chapter IV of \cite{MR0226684}) on $H$ if and only if $\varphi(a)=\exp(-\lambda(a))$, $a\in H$, with
\begin{equation}\label{101'}
\lambda(a)=-i\langle a,b\rangle+\frac{1}{2}\langle a,Ra\rangle-\int_{{H}}
\big(e^{i\langle a,x\rangle}-1-\frac{i\langle a,x\rangle}{1+\lVert x\rVert^2}\big)M(dx),
\end{equation}
where $b\in H$, $R\in\mathcal{L}_1^+(H)$ and $M$ is a L\'evy measure on $H$ (see e.g. Theorem 4.10, Chapter VI in \cite{MR0226684}), i.e. $M$ is a measure on $(H,\mathscr{B}(H))$ such that $M(\{0\})=0$ and $\int_H \big(1\wedge\rVert x\lVert^2\big)M(dx)<\infty$.\\

Now, let us recall the construction of a time-inhomogeneous Mehler semigroup by using its characteristic function. First we should state our hypotheses to be valid for the entire
paper:

\begin{description}
	\item ${\bf{(H1)}}$ $(U_{s,t})_{s\leq t}$ is a strong evolution family of uniformly bounded linear operators on $H$. 
	\item ${\bf{(H2)}}$ $\sigma:\mathbb{R}\rightarrow\mathcal{L}({H})$ is strongly continuous and bounded in operator norm.
	\item ${\bf{(H3)}}$ $\lambda:H\rightarrow\mathbb{C}$ is a negative definite and continuous function on $H$ with $\lambda(0)=0$ 
	and $\overline{\lambda(a)}=\lambda(-a)$ for all $a\in H$.
	\item ${\bf{(H4)'}}$ For all $s\leq t$
	\begin{equation*}
	a\longmapsto\exp{\bigg[-\int_{s}^{t}\lambda(\sigma^*(r)U_{r,t}^*a)dr\bigg]},\quad a\in H,
	\end{equation*}
	is Sazonov continuous, where $U^*$ denotes the adjoint of $U\in\mathcal{L}(H)$. 
\end{description}
Note that $(H3)$ does not imply the representation $(\ref{101'})$ for $\lambda$, unless we assume that $\lambda$ is Sazonov continuous. 

$(H1)-(H3)$ imply that, for all $s\leq t$, the function in $(H4)'$ is positive definite (see \cite{MR0481057}). Therefore, by the Minlos-Sazonov Theorem, they are characteristic functions of probability measures $\mu_{s,t}$ on $H$, i.e. we have
\begin{align}\label{4}
\widehat{\mu_{s,t}}(a)=\int_{H}e^{i\langle a,x\rangle}\mu_{s,t}(dx)=e^{-\int_s^t\lambda(\sigma^*(r)U_{r,t}^*a)dr},\quad a\in{H}.
\end{align}
If $(H1)-(H3)$ hold and $\lambda$ is itself Sazonov continuous, then $(H4)'$ holds automatically. This is easy to see as follows (see \cite{MR1745332}, \cite{MR2861314}). By $(\ref{101'})$ we have for all $a\in H$
\begin{equation*}
\begin{array}{ll}
\exp\bigg(-\int_s^t\lambda(\sigma^*(r)U_{r,t}^*a)dr\bigg)\\ \\
\quad\ =\exp\bigg\{\int_s^ti\langle a,U_{r,t}\sigma(r)b\rangle dr-\int_s^t\frac{1}{2}\langle \sigma^*(r)U_{r,t}^*a, R\sigma^*(r)U_{r,t}^*a\rangle dr\\
\quad\ +\int_s^t\int_H\big(e^{i\langle a,U_{r,t}\sigma(r)x\rangle}-1-\frac{i\langle a, U_{r,t}\sigma(r) x\rangle}{1+\lVert x\rVert^2}\big)M(dx)dr\bigg\} 
\end{array}
\end{equation*}
\begin{align}\label{R3}
\quad \quad \quad \quad\ \ =\exp\bigg(i\langle a,b_{s,t}\rangle-\frac{1}{2}\langle R_{s,t}a,a\rangle+\int_{{H}}\big(e^{i\langle a,x\rangle}-1-
\frac{i\langle a,x\rangle}{1+\lVert x\rVert^2}\big)M_{s,t}(dx)\bigg),
\end{align}
where 
\begin{align*}
R_{s,t}=\int_s^tU_{r,t}\sigma(r)R\sigma^*(r)U^*_{r,t}dr
\end{align*}
and
\begin{align*}
\begin{array}{ll}
b_{s,t}&=\int_s^tU_{r,t}\sigma(r)bdr\\ &+\int_s^t\int_{{H}}U_{r,t}\sigma(r)x\bigg(\frac{1}{1+\lVert U_{r,t}\sigma(r)x\rVert^2}-\frac{1}{1+\lVert x\rVert^2}\bigg)M(dx)dr
\end{array}
\end{align*}
are well-defined Bochner integrals with values in $\mathcal{L}_1^+(H)$ and $H$, respectively. In this formula, $M_{s,t}$ is a L\'evy measure on $H$, defined by:
\begin{equation}\label{R16}
M_{s,t}(B):=\int_s^tM\bigg((U_{r,t}\sigma(r))^{\ -1}(B\setminus\{0\})\bigg)dr,\quad B\in\mathscr{B}({H}).
\end{equation}
From representation $(\ref{R16})$ we immediately deduce by standard arguments that $(H4)'$ holds (see e.g. \cite{MR1745332}). However, as said before, we do not require that $\lambda$ is Sazonov, but we only assume $(H1)-(H3)$, $(H4)'$, resp. $(H4)$ below, in the entire paper.\\

Let $\pi_{s,t}(x,dy)$ be the translation of $\mu_{s,t}(dy)$ by $U_{s,t}x$, namely
\begin{equation}\label{110''}
\pi_{s,t}(x,dy)=\mu_{s,t}(dy-U_{s,t}x),\quad s\leq t,\ x\in H.
\end{equation}
We now show that the family $\pi=(\pi_{s,t})_{s\leq t}$ is a Markovian family of transition probabilities. By construction, all properties 
are obviously satisfied and only condition (iii) needs to be checked. By Proposition 2.2 in \cite{MR3466574}, (iii) is valid for $\pi$  if and only if 
\begin{equation}\label{00}
\mu_{s,t}=(\mu_{s,r}\circ U_{r,t}^{\ -1})*\mu_{r,t},\quad s\leq r\leq t,
\end{equation}
where $*$ is the convolution operator on $\mathscr{P}(H)$ (:=the set of all probability measures on $(H,\mathscr{B}(H))$). In terms of characteristic functions, $(\ref{00})$ is equivalent to: 
\begin{equation}\label{110'}
\widehat{\mu_{s,t}}(a)=\widehat{\mu_{s,r}}(U_{r,t}^{*}a)\widehat{\mu_{r,t}}(a),\quad a\in{H},\ s\leq r\leq t.
\end{equation}
But,
\begin{align*}
\widehat{\mu_{s,r}}(U_{r,t}^*a)\widehat{\mu_{r,t}}(a)&~=e^{-\int_s^r\lambda\big(\sigma^*(\ell)U_{\ell,r}^*(U_{r,t}^*a)\big)d\ell}\ e^{-\int_r^t\lambda
	\big(\sigma^*(\ell)U_{\ell,t}^*a\big)d\ell}\\
&~=e^{-\int_s^r\lambda\big(\sigma^*(\ell)(U_{r,t}U_{\ell,r})^*a\big)d\ell}\ e^{-\int_r^{t}\lambda\big(\sigma^*(\ell)U_{\ell,t}^*a\big)d\ell}\\
&~=e^{-\int_s^r\lambda\big(\sigma^*(\ell)U_{\ell,t}^*a\big)d\ell}\ e^{-\int_r^{t}\lambda\big(\sigma^*(\ell)U_{\ell,t}^*a\big)d\ell}\\
&~=e^{-\int_s^{t}\lambda\big(\sigma^*(\ell)U_{\ell,t}^*a\big)d\ell}\\
&~=\widehat{\mu_{s,t}}(a).
\end{align*}
Hence, $\pi=(\pi_{s,t})_{s\leq t}$ is a Markovian family of transition probabilities and as we mentioned before, it is called  
\emph{the time-inhomogeneous generalized Mehler semigroup}.
The characteristic function of $\pi_{s,t}(x,dy)$ is for $x\in H$ given by
\begin{equation}\label{11!}
\int e^{i\langle a,y\rangle}\pi_{s,t}(x,dy)=e^{i\langle a,U_{s,t}x\rangle-\int_s^t\lambda(\sigma^*(r)U_{r,t}^*a)dr},\quad a\in H.
\end{equation}
\begin{remark}\label{R2.1}
Conditions $(H3)$ and $(H4)'$ are crucial to ensure that the probability measures $\mu_{s,t}$, $s\leq t$, exist, hence that the time-inhomogeneous generalized Mehler semigroup in \eqref{110''} exists.
\end{remark}
\section{Extremal entrance laws}
Let $\nu=(\nu_t)_{t\in\mathbb{R}}\subset\mathscr{P}(H)$ be such that for all $t\in\mathbb{R}$, $a\in H$  
\begin{equation}\label{1}
\int_{H}\lvert\langle a,x\rangle\rvert\ \nu_t(dx)<\infty,\quad t\in\mathbb{R}.
\end{equation}
Since for each $t\in\mathbb{R}$, $\nu_t$ is a probability, hence a finite measure, the uniform boundedness principle implies that the linear functional
$a\mapsto\int_{H}\langle a,x\rangle\nu_t(dx)$ is continuous on ${H}$.
Hence, by the Riesz representation theorem there exists $\kappa_t\in H$ such that 
\begin{equation}\label{2}
\int_{H}\langle a,x\rangle \nu_t(dx)=\langle a,\kappa_t\rangle,\quad a\in{H},\ t\in\mathbb{R},
\end{equation}
i.e. $\kappa_t$ is the mean of $\nu_t$.\\
We recall 
\begin{align*}
\mathcal{K}(\pi):=\bigg\{\nu:=(\nu_t)_{t\in\mathbb{R}}\in\mathscr{P}({H})^\mathbb{R}\ \bigg\lvert\ 
\int_H\pi_{s,t}(x,B)\nu_s(dx)=\nu_{t}(B),&~\ s\leq t,\ s,t\in\mathbb{R}\\ &~\ ,B\in\mathscr{B}(H)\bigg\}
\end{align*}
and 
\begin{align*}
\mathcal{K}(U):=\big\{\kappa=(\kappa_t)_{t\in\mathbb{R}}\in H^\mathbb{R}\mid U_{s,t}\kappa_s=\kappa_{t},\ s\leq t,\ s,t\in\mathbb{R}\big\}.
\end{align*}
\begin{remark}\label{R5}
	{\emph{(i)}} Clearly, by the strong continuity of $U_{s,t}$, $s\leq t$, we have $\mathcal{K}(U)\subset C(\mathbb{R};H)$.\\
	{\emph{(ii)}} Let $\nu=(\nu_t)_{t\in\mathbb{R}}\in\mathcal{K}(\pi)$ and $a\in H$ such that there exist $t_n\in\mathbb{R}$, $n\in\mathbb{N}$, $t_{n+1}\leq t_n$, $\lim_{n\rightarrow\infty}t_n=-\infty$ and 
	\begin{align*}
	\int_{H}\lvert\langle a,x\rangle\rvert\nu_{t_n}(dx)<\infty,\quad \forall n\in\mathbb{N}.
	\end{align*}
	Then obviously $\nu$ satisfies $(\ref{1})$.\\
	{\emph{(iii)}} Obviously, 
	\begin{align*}
	(\nu_t)_{t\in\mathbb{R}}\in\mathcal{K}(\pi)\Longleftrightarrow \widehat{\nu_t}(a)=\widehat{\nu_s}(U^*_{s,t}a)\ \widehat{\mu_{s,t}}(a),\quad \forall s\leq t,\ a\in H.
	\end{align*} 	
\end{remark}

We also recall that $\mathcal{K}^1(\pi)$ is the set of all $\nu=(\nu_t)_{t\in\mathbb{R}}\in\mathcal{K}(\pi)$ which have finite weak first moments, i.e. satisfy $(\ref{1})$ for all $t\in\mathbb{R}$, $a\in H$. \\
The map $\nu\rightarrow\kappa$ from $\mathcal{K}^1(\pi)$ to $H^{\mathbb{R}}$ is denoted by $\bf{p}$. This $\bf{p}(\nu)$ is just the mean of $\nu$ and 
is called the {\it{projection of $\nu$}} in \cite{MR1013934}.\\

Note that $(H1)-(H3),\ (H4)'$ are still in force. In addition, from now on we also assume the following hypotheses: 
\begin{description}
	\item ${\bf{(H4)}}$ For all $t\in\mathbb{R}$, $r\longmapsto\lambda(\sigma^*(r)U_{r,t}^*a)$ is Lebesgue integrable on $(-\infty,t)$ for all $a\in H$ and
	\begin{equation}\label{R15}
	a\longmapsto\exp{\bigg[-\int_{-\infty}^{t}\lambda(\sigma^*(r)U_{r,t}^*a)dr\bigg]},\quad a\in H,
	\end{equation}
	is Sazonov continuous for all $t\in\mathbb{R}$. Furthermore, the probability measure $\mu_{-\infty,t}$ defined by
	\begin{align*}
	\widehat{\mu_{-\infty,t}}(a):=e^{-\int_{-\infty}^t\lambda(\sigma^*(r)U_{r,t}^*a)dr},\quad a\in H,
	\end{align*}
	has finite weak first moments for all $t\in\mathbb{R}$.
	\item ${\bf{(H5)}}$ $\lambda=\overline{\lambda}$.
\end{description}

\begin{remark}\label{R18}
	{\emph{(i)}} Obviously, $(H4)$ is stronger than $(H4)'$. Conditions $(H4)$ and $(H5)$ are cruical to ensure that $\mathcal{K}^1(\pi)$ is not empty (see Lemma \ref{R4} and its proof below). Since $\lambda$ itself is not assumed to be Sazonov continuous, cases with cylindrical Levy processes in \eqref{2'} (in particular, cylindrical Wiener processes with e.g.~$\lambda(a)=\|a\|^2$, $a\in H$) are covered provided $U_{s,t}$, $s\leq t$, are Hilbert-Schmidt (see \cite[Section 8]{MR2032114} for an example, namely  the stochatic heat equation with cylindrical Levy noise).\\
	{\emph{(ii)}} Suppose that $\lambda$ from $(H3)$ is Sazonov continuous, or equivalently $\lambda$ has a representation as in $(\ref{101'})$. Then obviously $(H5)$ holds if and only if
	\begin{align}\label{R14}
	\lambda(a)=\frac{1}{2}\langle a,Ra\rangle+\int_H(1-\cos\langle a,x\rangle)M(dx),
	\end{align}
	for all $a\in H$. Furthermore, $M$ is symmetric in this case.\\
	{\emph{(iii)}} Suppose that $(H1)-(H3)$, $(H5)$ hold. Now we formulate additional (checkable) assumptions on $\lambda$ from $(H3)$ and $(U_{s,t})_{s\leq t}$ from $(H2)$, which imply that $(H4)$ holds. So, about $\lambda$ we additionally assume:
	\begin{description}
		\item $(\lambda.1)$ $\lambda$ is Sazonov continuous, or equivalently $\lambda$ has a representation as in $(\ref{R14})$ with corresponding L\'evy measure $M$.
		\item $(\lambda.2)$ 
		\begin{align*}
		\int_{\{\lVert\cdot\rVert>1\}}\lVert x\rVert\ M(dx)<\infty.
		\end{align*}
	\end{description}
	Furthermore, assume on $(U_{s,t})_{s\leq t}$ from $(H1)$:
	\begin{description}
		\item $(U.1)$ There exist $c,\omega\in(0,\infty)$ such that 
		\begin{align*}
		\lVert U_{s,t}\rVert_{\mathcal{L}(H)}\leq ce^{-\omega(t-s)},\quad \forall s\leq t.
		\end{align*} 
	\end{description}
	Then by the same arguments as those implying the representation $(\ref{R3})$, hence the Sazonov continuity  of the function in $(\ref{R3})$ with $-\infty$ replacing $s$, show that the function in $(\ref{R15})$ has the representation $(\ref{R3})$ with $-\infty$ replacing $s$, and hence is Sazonov continuous.\\ The only difference is that, we need to check that $b_{-\infty,t}$ is well-defined. This, however, immediately follows from conditions $(U.1)$ and $(\lambda.2)$, since for all $r\leq t$, $x\in H$ and $C_{\sigma}:=c\ \sup_{r\in\mathbb{R}}\lVert \sigma(r)\rVert_{\mathcal{L}(H)}$
	\begin{align*}
	&~\lVert U_{r,t}\sigma(r) x\rVert\ \frac{\lVert x\rVert^2+\lVert U_{r,t}\sigma(r) x\rVert^2}{(1+\lVert U_{r,t}\sigma(r)x\rVert^2)(1+\lVert x\rVert^2)}\\ \\
	\leq&~\lVert U_{r,t}\sigma(r) x\rVert\ (1+e^{-2\omega(t-r)}C_\sigma^2)\ \frac{\lVert x\rVert^2}{1+\lVert x\rVert^2}\\ \\
	\leq&~\lVert x\rVert\ e^{-\omega(t-r)}C_\sigma(1+C_\sigma^2)\ (\lVert x\rVert^2 \wedge 1).
	\end{align*}
	Hence by the Minlos-Sazonov Theorem, the measures $\mu_{-\infty,t}$, $t\in\mathbb{R}$, in $(H4)$ exist.
	
	To obtain that $(H4)$ holds, it remains to show each $\mu_{-\infty,t}$ has finite weak first moments. To show this, it suffixes to consider the case $R=0$, because if not, we just have to convolute with $\mathcal{N}(0,R)$, i.e. the centered Gaussian measure with covariance operator $R$, which has all strong moments, so the convolution, in particular, will preserve finite weak first moments.
	
	Since for $a\in H$ the L\'evy measure of $\mu_{-\infty,t}\circ\langle a,\cdot\rangle^{-1}$ is $M_{-\infty,t}\circ\langle a,\cdot\rangle^{-1}$ (with $M_{-\infty,t}$ defined as in $(\ref{R16})$ with $s=-\infty$), it follows by conditions $(\lambda.2)$ and $(U.1)$ that each $\mu_{-\infty,t}$ has finite first weak moments, so $(H4)$ holds. Indeed, by \cite{MR1739520}, Theorem 25.3, we only need to check that 
	\begin{align*}
	\int_{\{\lvert\cdot\rvert>1\}}\lvert s\rvert\ \big(M_{-\infty,t}\circ\langle a,\cdot\rangle^{-1}\big)(ds)<\infty.
	\end{align*}
	But by the definition of $M_{-\infty,t}$, the left hand side is equal to 
	\begin{align*}
	&~\int_{-\infty}^{t}\int_{\{\lvert\langle\sigma^{*}(r)U_{r,t}^{*}a,\cdot\rangle\rvert>1\}}\big\lvert\langle\sigma^{*}(r)U_{r,t}^{*}a,x\rangle\big\rvert M(dx)dr\\ \\
	\leq&~\int_{-\infty}^{t}\int_{\{\lVert\cdot\rVert\leq 1\}}\langle\sigma^{*}(r)U_{r,t}^{*}a,x\rangle^2 M(dx)dr\\
	&~+\int_{-\infty}^{t}\int_{\{\lVert\cdot\rVert>1\}}C_\sigma e^{-\omega(t-r)}\lVert a\rVert\lVert x\rVert M(dx)dr\\ \\
	\leq&~\frac{1}{2\omega}\ C_\sigma^2\ \lVert a\rVert^2\int_{\{\lVert\cdot\rVert\leq1\}}\lVert x\rVert^2 M(dx)+ 
	\frac{1}{\omega}\ C_\sigma\ \lVert a\rVert\int_{\{\lVert\cdot\rVert>1\}}\lVert x\rVert M(dx),
	\end{align*}
	which is finite by $(\lambda.2)$.
\end{remark}
In Section 5, we shall give explicit examples for $\lambda$ satisfying $(H3)$, $(H5)$, $(\lambda.1)$ and $(\lambda.2)$, hence $(H4)$.
\begin{lemma}\label{R4}
	$(\mu_{-\infty,t})_{t\in\mathbb{R}}\in\mathcal{K}^1(\pi)$ with $\kappa_t=0$ for all $t\in\mathbb{R}$.
	\begin{proof}
		Analogous to the proof of $(\ref{110'})$, and $\kappa_t=0$ for all $t\in\mathbb{R}$, follows, since by $(H5)$ the Fourier transform $\widehat{\mu_{-\infty,t}}$ is real, hence $\mu_{-\infty,t}$ is symmetric (i.e. $\mu_{-\infty,t}(dx)=\mu_{-\infty,t}(-dx)$ for all $t\in\mathbb{R}$).
	\end{proof}
\end{lemma}

\begin{lemma}\label{49'}
	We have for all $s\leq t$ and $a\in H$:
	\begin{description}
		\item {\emph{(i)}}  
		\begin{align*}
		\int_H\lvert \langle a,y\rangle\rvert\ \mu_{s,t}(dy)<\infty\ \ \text{and}\ \ \int_H \langle a,y\rangle\ \mu_{s,t}(dy)=0.
		\end{align*}
		\item {\emph{(ii)}} 
		\begin{align}\label{12!}
		\int_H\lvert\langle a,y\rangle\rvert\ \pi_{s,t}(x,dy)<\infty\ \ \text{and}\ \ \int_H\langle a,y\rangle\ \pi_{s,t}(x,dy)=\langle a, U_{s,t}x\rangle
		\end{align}
		for all $x\in H$.
	\end{description}
	\begin{proof} {\emph{(i)}}: For all $x\in H$, we have
		\begin{align*}
		\int_{H}\lvert\langle a,y\rangle\rvert\ \mu_{s,t}(dy)
		\leq\int_{H}\lvert\langle a,y+U_{s,t}x\rangle\rvert\ \mu_{s,t}(dy)+\int_{H}\lvert\langle a,U_{s,t}x\rangle\rvert\ \mu_{s,t}(dy).
		\end{align*}
		By integrating over $x$ with respect to $\mu_{-\infty,s}$ and using Lemma \ref{R4} as well as $(H4)'$, we get
		\begin{align*}
		\int_{H}\lvert\langle a,y\rangle\rvert\ \mu_{s,t}(dy)&~\leq\int_{H}\int_H\lvert\langle a,y\rangle\rvert\ \pi_{s,t}(x,dy)\ \mu_{-\infty,s}(dx)+\int_H\lvert\langle a,U_{s,t}x\rangle\rvert\ \mu_{-\infty,s}(dx)\\
		&~=\int_H\lvert\langle a,y\rangle\rvert\ \mu_{-\infty,t}(dy)+\int_H\lvert\langle U^*_{s,t}a,x\rangle\rvert\ \mu_{-\infty,s}(dx)<\infty
		\end{align*}
		Thus, {\emph{(i)}} holds, because each $\mu_{s,t}$ is symmetric.\\
		{\emph{(ii)}} immediately follows from {\emph{(i)}}.
	\end{proof}
\end{lemma}

\begin{proposition}\label{33!}
	Assume $(H1)-(H5)$. Then for each $\nu\in\mathcal{K}^1(\pi)$, $\kappa:={\bf{p}}(\nu)\in\mathcal{K}(U)$.
	\begin{proof}
		Let $a\in H$. We need to check that $\int_{{H}}\langle a,x\rangle \nu_t(dx)=\langle a,U_{s,t}\kappa_s\rangle$ for all $s\leq t$. By the definition of 
		$\mathcal{K}(\pi)$, we get
		\begin{align*}
		\int_{{H}}\langle a,x\rangle \nu_t(dx)&~=\int_{{H}}\bigg(\int_{{H}}\langle a,y\rangle\pi_{s,t}(x,dy)\bigg)\nu_s(dx).
		\end{align*}
		Lemma \ref{49'} implies
		\begin{align*}
		\int_{{H}}\bigg(\int_{{H}}\langle a,y\rangle\pi_{s,t}(x,dy)\bigg)\nu_s(dx)
		&~=\int_{{H}}\langle a,U_{s,t}x\rangle \nu_s(dx)\\
		&~=\int_{{H}}\langle U_{s,t}^*a,x\rangle \nu_s(dx)\\
		&~=\langle U_{s,t}^*a,\kappa_s\rangle=\langle a,U_{s,t}\kappa_s\rangle,
		\end{align*}
		which completes the proof.
	\end{proof}
\end{proposition}
As a part of our main result (see Theorem \ref{27!} below) we shall obtain that $\mathcal{K}^1(\pi)$ is a simplex, i.e. that each element in  $\mathcal{K}^1(\pi)$ has a unique representation as an integral over its extreme points  $\mathcal{K}^1_e(\pi)$.
The next result is a first step in this direction and in its proof we also identify the difficulty why this is not a trivial consequence of E. B. Dynkin's result in \cite{MR0518321}, which as recalled in the introduction, states that  $\mathcal{K}(\pi)$ is a simplex.
\begin{proposition}\label{R8}
	{\emph{(i)}}  $\mathcal{K}^1_e(\pi)\subset\mathcal{K}_e(\pi)$.\\
	{\emph{(ii)}} Let $A\subset H$ be a countable $\mathbb{Q}$-vector space such that $A$ is dense in $H$ (in the norm topology). Let $H_0:=span\ A$ be its $\mathbb{R}$-linear span. Define
	\begin{align*}
	\mathcal{K}_e^{H_0}(\pi):=\bigg\{\nu=(\nu_t)_{t\in\mathbb{R}}\in\mathcal{K}_e(\pi)\bigg\lvert\int_H\lvert\langle a,x\rangle\rvert\ \nu_t(dx)<\infty,\quad \forall t\in\mathbb{R},\ a\in H_0\bigg\}.
	\end{align*}
	Then\\
	{\emph{(a)}} $\mathcal{K}_e^{H_0}(\pi)=\big\{\nu=(\nu_t)_{t\in\mathbb{R}}\in\mathcal{K}_e(\pi)\big\lvert\int_H\lvert\langle a,x\rangle\rvert\ \nu_{-n}(dx)<\infty,\quad \forall n\in\mathbb{N},\ a\in A\big\}$.\\
	{\emph{(b)}} Let $\nu=(\nu_t)_{t\in\mathbb{R}}\in\mathcal{K}^1(\pi)$. Then $\nu$ has a unique representation as an integral 
	\begin{align*}
	\nu=\int_{\mathcal{K}_e^{H_0}(\pi)}\widetilde{\nu}\ \xi_\nu(d\widetilde{\nu})
	\end{align*}
	over $\mathcal{K}_e^{H_0}(\pi)$.
	\begin{proof}
		{\emph{(i)}}: Let $\nu\in\mathcal{K}_e^1(\pi)$ and $\nu^{(1)},\nu^{(2)}\in\mathcal{K}(\pi)$ with $\nu^{(1)}\neq\nu^{(2)}$ such that  $\nu=\alpha\nu^{(1)}+(1-\alpha)\nu^{(2)}$
		for some 
		$\alpha\in(0,1)$. Then $\nu^{(1)}\leq \frac{1}{\alpha}\nu$ and $\nu^{(2)}\leq \frac{1}{1-\alpha}\nu$. Hence $\nu^{(1)},\nu^{(2)}\in\mathcal{K}^1(\pi)$. Therefore, $\nu^{(1)}=\nu^{(2)}$, which means that $\nu\in\mathcal{K}_e(\pi)$.  \\
		{\emph{(ii)}}: {\emph{(a)}} follows by linearity and Remark \ref{R5} {\emph{(ii)}}, so let us prove {\emph{(b)}}. As mentioned before, by \cite{MR0518321}, $\mathcal{K}(\pi)$ is a simplex, so each element in $\mathcal{K}(\pi)$ has a unique representation as an integral over its extreme points $\mathcal{K}_e(\pi)$. More precisely, consider the $\sigma$-algebra $\mathscr{A}$ on $\mathcal{K}_e(\pi)$, generated by all maps 
		\begin{align*}
		\mathcal{K}_e(\pi)\ni(\nu_s)_{s\in\mathbb{R}}\ \mapsto\ \nu_t\in\mathscr{P}(H),\quad t\in\mathbb{R},
		\end{align*}
		where $\mathscr{P}(H)$ is equipped with the $\sigma$-algebra generated by the weak topology. Then for each $\nu=(\nu_t)_{t\in\mathbb{R}}\in\mathcal{K}(\pi)$, there exists a unique probability measure $\xi_\nu$ on $(\mathcal{K}_e(\pi),\mathscr{A})$ such that 
		\begin{align}\label{R6}
		\nu=\int_{\mathcal{K}_e(\pi)}\widetilde{\nu}\ \xi_\nu(d\widetilde{\nu}).
		\end{align}
		Let $\nu\in\mathcal{K}^1(\pi)$. Then for all $t\in\mathbb{R}$, $a\in H$
		\begin{align*}
		\infty>\int_H\lvert\langle a,x\rangle\rvert\ \nu_t(dx)=\int_{\mathcal{K}_e(\pi)}\int_H\lvert\langle a,x\rangle\rvert\ \widetilde{\nu_t}(dx)\  \xi_\nu(d\widetilde{\nu}),
		\end{align*}
		which yields 
		\begin{align*}
		\int_H\lvert\langle a,x\rangle\rvert\ \widetilde{\nu_t}(dx)<\infty,
		\end{align*}
		for $\xi_\nu-a.e.$ $\widetilde{\nu}\in\mathcal{K}_e(\pi)$. Here initially the $\xi_\nu$-zero set depends on $t$ and $a$. But specializing to $t=-n$, $n\in\mathbb{N}$, and $a\in H_0$ by assertion {\emph{(ii)}} part {\emph{(a)}} it can be chosen independent of $t\in\mathbb{R}$ and $a\in H_0$. Hence
		\begin{align*}
		\xi_\nu\big(\mathcal{K}_e^{H_0}(\pi)\big)=1
		\end{align*}
		and $(\ref{R6})$ holds with $\mathcal{K}_e^{H_0}(\pi)$ replacing $\mathcal{K}_e(\pi)$, which is the assertion, since the uniqueness is obvious from the uniqueness of the representation of $\nu$ over $\mathcal{K}_e(\pi)$. So, {\emph{(ii)}} part {\emph{(b)}} is proved.
	\end{proof}
\end{proposition}

Now we are able to prove the main result of this paper. Before, we need to define the Markov processes associated with $\pi$. \\
For a given $\nu\in\mathcal{K}(\pi)$, one can construct a unique probability measure $\mathbb{P}_\nu$ on the space $\varOmega:=H^{\mathbb{R}}$ with $\sigma$-algebra $\mathcal{F}:=\sigma(X_t\mid t\in\mathbb{R})$ such that for all $t_1\leq \cdots \leq t_n$
\begin{align*}
\mathbb{P}_\nu\big[X_{t_1}\in dx_1,,..., X_{t_n}\in dx_n\big]:=\pi_{t_{n-1},t_{n}}(x_{n-1},dx_n)\cdots\pi_{t_1,t_2}(x_1,dx_2)\nu_{t_1}(dx_1),
\end{align*}
where $X_t:\varOmega\rightarrow H$ is the canonical coordinate process. Obviously, $\nu\mapsto\mathbb{P}_\nu$ is then convex and injective, since 
\begin{align*}
\mathbb{P}_\nu\circ X_t^{-1} = \nu_t,\quad t\in\mathbb{R},
\end{align*}
i.e. $\nu_t$, $t\in\mathbb{R}$, are the one dimensional marginals of $\mathbb{P}_\nu$. Furthermore, this $\mathbb{P}_\nu$ is Markovian, i.e.,
\begin{align}\label{m}
\mathbb{P}_\nu[X_t\in dz\mid \mathcal{F}_s] = \pi_{s,t}(X_s, dz)\quad \forall t, s\in \mathbb{R},\ t > s,
\end{align}
where $\mathcal{F}_s:=\sigma(X_r\mid r\leq s)$. Define the convex set $\mathcal{M}(\pi):=\{\mathbb{P}_\nu\mid\nu\in\mathcal{K}(\pi)\}$.

\begin{lemma}\label{pp}
	Let $\mathbb{P}\in\mathcal{M}(\pi)$. Then $\mathbb{P}$ is an extremal point of $\mathcal{M}(\pi)$ if and only if $\mathbb{P}(\varGamma)=1$ or 0 for every $\varGamma\in\mathcal{F}_{-\infty}:=\bigcap_{s\in\mathbb{R}}\mathcal{F}_s$.
	\begin{proof}
		See the proof of Lemma 2.4. in \cite{MR1159586}.
	\end{proof}
\end{lemma}

\begin{theorem}\label{27!}
	Let $(\pi_{s,t})_{s\leq t}$ be the time-inhomogeneous (generalized) Mehler semigroup on $H$ as above. Assume that $(H1)-(H5)$ hold. \\
	{\emph{a)}} Let $\kappa=(\kappa_t)_{t\in\mathbb{R}}\in\mathcal{K}(U)$. Then
	\begin{align*}
	\nu^\kappa(dy):=\big(\mu_{-\infty,t}(dy-\kappa_t)\big)_{t\in\mathbb{R}}\in\mathcal{K}^1_e(\pi),
	\end{align*}
	and ${\bf{p}}(\nu^\kappa)=\kappa$. Here, $\mu_{-\infty,t}(dy-\kappa_t)$ denotes the image measure of $\mu_{-\infty,t}$ under the map $H\ni x\mapsto x+\kappa_t$, $t\in\mathbb{R}$.\\
	{\emph{b)}} The map
	\begin{align*}
	\mathcal{K}(U)\ni\kappa=(\kappa_t)_{t\in\mathbb{R}}\ \longmapsto\ \nu^\kappa\in\mathcal{K}^1_e(\pi)
	\end{align*}
	is a bijection.\\
	{\emph{c)}} $\mathcal{K}^1(\pi)$ is a simplex, i.e. each $\nu\in\mathcal{K}^1(\pi)$ has a unique representation as an integral 
	\begin{align*}
	\nu=\int_{\mathcal{K}_e^1(\pi)}\widetilde{\nu}\ \xi_\nu(d\widetilde{\nu})
	\end{align*}
	over its extreme points $\mathcal{K}_e^1(\pi)$.
	\begin{proof} The following claims {\emph{(i)}}, {\emph{(ii)}} and {\emph{(iii)}} together with Proposition \ref{R8} prove the theorem.
		\begin{description}
			\item {\it{Claim (i)}} $\nu^\kappa\in\mathcal{K}^1(\pi)$.
		\end{description}
		{\it{Proof}}: Let $\kappa\in\mathcal{K}(\pi)$, $t\in\mathbb{R}$. Then 
		\begin{align}\label{09!}
		\int_{H}e^{i\langle a,y\rangle}\nu^\kappa_t(dy)=e^{i\langle a,\kappa_t\rangle-\int_{-\infty}^t\lambda(\sigma^*(r)U_{r,t}^*a)dr},\quad \forall t\in\mathbb{R},\ a\in H.
		\end{align}
		Since $\mu_{-\infty,t}$ has finite weak first moments, so has $\nu_t^\kappa$. It remains to prove that $(\nu_t^\kappa)_{t\in\mathbb{R}}$ belongs to $\mathcal{K}(\pi)$. But for all $a\in H$
		\begin{align*}
		\widehat{(\nu_s^\kappa \pi_{s,t}\ )}(a)=&~\int_H\bigg(\int_He^{i\langle a,y\rangle}\pi_{s,t}(x,dy)\bigg)\nu^\kappa_s(dx)\\ 
		=&~\int_{{H}}e^{i\langle a,U_{s,t}x\rangle-\int_s^t\lambda\big(\sigma^*(r)U_{r,t}^*a\big)dr}\nu^\kappa_s(dx)\\ 
		=&~\int_{{H}}e^{i\langle U_{s,t}^*a,x\rangle}\nu^\kappa_s(dx).\ e^{-\int_s^t\lambda\big(\sigma^*(r)U_{r,t}^*a\big)dr}\\ 
		=&~e^{i\langle U_{s,t}^*a,\kappa_s\rangle-\int_{-\infty}^s\lambda\big(\sigma^*(r)U_{r,s}^*(U_{s,t}^*a)\big)dr}.\ e^{-\int_s^t\lambda\big(\sigma^*(r)U_{r,t}^*a\big)dr}\\ 
		=&~e^{i\langle a,U_{s,t}\kappa_s\rangle-\int_{-\infty}^s\lambda\big(\sigma^*(r)\ (U_{s,t}U_{r,s})^*a\big)dr-\int_s^t\lambda\big(\sigma^*(r)U_{r,t}^*a\big)dr}\\ 
		=&~e^{i\langle a,\kappa_{t}\rangle-\int_{-\infty}^s\lambda\big(\sigma^*(r)U_{r,t}^*a\big)dr-\int_s^{t}\lambda\big(\sigma^*(r)U_{r,t}^*a\big)dr}\\ 
		=&~e^{i\langle a,\kappa_{t}\rangle-\int_{-\infty}^{t}\lambda\big(\sigma^*(r)U_{r,t}^*a\big)dr}\\ 
		=&~\widehat{\nu^\kappa_{t}}(a).
		\end{align*}
		Hence Claim {\emph{(i)}} is proved.
		\begin{description}
			\item {\it{Claim (ii)}} {\emph{Let $\mathcal{K}_e^{H_0}(\pi)$ be defined as in Proposition \ref{R8} and let $\nu\in\mathcal{K}^{H_0}_e(\pi)$. Then $\nu\in\mathcal{K}^1_e(\pi)$. Define $\kappa:=\bf{p}(\nu)$ ($\in\mathcal{K}(U)$ by Proposition \ref{33!}). Then 
					\begin{align}\label{R26}
					\widehat{\nu_t}(a)=e^{i\langle a,\kappa_t\rangle-\int_{-\infty}^t\lambda(\sigma^*(r)U_{r,t}^*a)dr},\quad \forall a\in H,\ t\in\mathbb{R}.
					\end{align}
					In particular, ${\bf{p}}\mid_{\mathcal{K}_e^1(\pi)}:\ \mathcal{K}_e^1(\pi)\rightarrow \mathcal{K}(U)$ is injective.}}
		\end{description}
		{\it{Proof}}: Since injective convex mappings map extreme points to extreme points, Lemma \ref{pp} implies that $\mathbb{P}_\nu$ is trivial on $\mathcal{F}_{-\infty}$. Thus, for every $t\in\mathbb{R}$ and every measurable function $f$ on ${H}$ with $\mathbb{E}_\nu \lvert f(X_t)\rvert<\infty$, 
		we have
		\begin{align}\label{07!}
		\begin{array}{lll}
		\int_H f(y)\nu_t(dy)&~=\mathbb{E}_\nu f(X_t)\\ \\
		&~=\mathbb{E}_\nu\big[f(X_t)\arrowvert\mathcal{F}_{-\infty}\big]\\ \\
		&~=\lim_{n\rightarrow \infty} \mathbb{E}_\nu\big[f(X_t)\arrowvert\mathcal{F}_{-n}\big]\\ \\
		&~=\lim_{n\rightarrow \infty} \int_H f(y)\pi_{-n,t}(X_{-n},dy),\quad \mathbb{P}_\nu-a.s..
		\end{array}
		\end{align}
		Note that in the third line we applied the backwards Martingale convergence theorem to the process $\mathbb{E}_\nu\{f(X_t)\arrowvert\mathcal{F}_{-n}\}$, $-n\leq t$, which is a martingale. Furthermore, in the fourth line we used the Markov property $(\ref{m})$ of our process.\\
		Now, $(\ref{07!})$ and $(\ref{11!})$, imply for every $a\in H$ that $\mathbb{P}_\nu-a.s.$ 
		\begin{align}\label{R9}
		\widehat{\nu_t}(a)&~=\lim_{n\rightarrow \infty}\int_H e^{i\langle a,y\rangle}\pi_{-n,t}(X_{-n},dy)\nonumber\\ 
		&~=\lim_{n\rightarrow \infty}e^{i\langle U_{-n,t}^*a,X_{-n}\rangle-\int_{-n}^{t}\lambda(\sigma^*(r)U_{r,t}^*a)dr}\nonumber\\
		&~=\lim_{n\rightarrow \infty}e^{i\langle U_{-n,t}^*a,X_{-n}\rangle}\ e^{-\int_{-\infty}^{t}\lambda(\sigma^*(r)U_{r,t}^*a)dr}.
		\end{align}
		And finally, by applying $(\ref{12!})$ and  $(\ref{07!})$, because $\nu\in\mathcal{K}_e^{H_0}(\pi)$, we obtain that for all $t\in\mathbb{R}$, $a\in H_0$
		\begin{align*}
		&~\langle U_{-n,t}^*a,X_{-n}\rangle=\int_H \langle a,y\rangle\pi_{-n,t}(X_{-n},dy)\xrightarrow{n\rightarrow\infty}
		\int_H \langle a,y\rangle\nu_t(dy)=
		\langle a,\kappa_t\rangle,\quad \mathbb{P}_\nu-a.s..
		\end{align*}
		This and $(\ref{R9})$ imply that $\forall t\in\mathbb{R},\ a\in H_0$
		\begin{align}\label{R10}
		\widehat{\nu_t}(a)=e^{i\int_H\langle a,y\rangle\nu_t(dy)}\ e^{-\int_{-\infty}^{t}\lambda(\sigma^*(r)U_{r,t}^*a)dr}.
		\end{align}
		We now show that $(\ref{R9})$ and $(\ref{R10})$ imply that $\nu\in\mathcal{K}_e^1(\pi)$ and that $(\ref{R26})$ holds. So, fix $t\in\mathbb{R},\ n\in\mathbb{N}$, and let $\{e_i\mid i\in\mathbb{N}\}\subset H_0$ be an orthonormal basis of $H$. Define $P_n:H\rightarrow H_n:=span\{e_1,\cdots,e_n\}$ and 
		\begin{align*}
		\mathring{\nu}_t^n:=\nu_t\circ P_n^{-1},\quad \mathring{\mu}^n_{-\infty,t}:=\mu_{-\infty,t}\circ P_n^{-1}.
		\end{align*}
		We extend these measures on $\mathscr{B}(H_n)$ by zero to $\mathscr{B}(H)$, i.e. we define for $B\in\mathscr{B}(H)$
		\begin{align*}
		\nu_t^n(B)&~:=\mathring{\nu}_t^n(B\cap H_n),\\
		\mu_{-\infty,t}^n(B)&~:=\mathring{\mu}_{-\infty,t}^n(B\cap H_n).
		\end{align*}
		Note that then for any $f:H\rightarrow\mathbb{C}$ bounded, $\mathscr{B}(H)$-measurable
		\begin{align*}
		\int_H f\ d\nu_t^n = \int_H f\circ P_n\ d\nu_t
		\end{align*}
		and likewise for $\mu_{-\infty,t}^n$. Then, in particular, we have for all $a\in H$
		\begin{align*}
		\widehat{\nu^n_t}(a)=\widehat{\nu_t}(P_n a),\quad \widehat{\mu^n_{-\infty,t}}(a)=\widehat{\mu_{-\infty,t}}(P_n a).
		\end{align*}
		Thus $(\ref{R10})$ implies 
		\begin{align}\label{R11}
		\widehat{\nu^n_t}(a)=e^{i\langle a,\kappa_t^n\rangle}\ \widehat{\mu^n_{-\infty,t}}(a),\quad \forall a\in H,
		\end{align}
		where 
		\begin{align*}
		\kappa_t^n:=\int_H y\ \nu_t^n(dy)=\int_{H_n} y\ \mathring{\nu}_t^n(dy)\in H_n\subset H.
		\end{align*}
		Letting $n\rightarrow\infty$ in $(\ref{R11})$, we obtain that for all $a\in H$
		\begin{align*}
		F(a):=\lim_{n\rightarrow\infty}e^{i\langle a,\kappa_t^n\rangle}=\frac{\widehat{\nu_t}(a)}{\widehat{\mu_{-\infty,t}}(a)}
		\end{align*}
		exists, is positive definite and Sazonov continuous with $F(0)=0$. Hence, by the Minlos-Sazonov Theorem, there exists a probability measure $\mu$ on $\mathscr{B}(H)$ such that 
		\begin{align*}
		\widehat{\mu}(a)=F(a),\quad a\in H,
		\end{align*}
		and, thus, by \cite{MR1435288}, Chap. IV, Proposition 3.3, the sequence of Dirac measures $\delta_{\kappa_t^n}$, $n\in\mathbb{N}$, converges weakly to $\mu$ with respect to the weak topology on $H$. From this, it is easy to show that, there exists $\kappa_t\in H$ such that $\kappa^n_t\rightharpoonup\kappa_t$ (i.e. weakly) in $H$ as $n\rightarrow\infty$. Indeed, for $a\in H$, let $\chi\in C_b(\mathbb{R})$, $\chi_N=1$ on $[-N,N]$, $\chi_N=0$ on $\mathbb{R}\setminus(-(N+1),N+1)$. Then
		\begin{align*}
		\lim_{n\rightarrow\infty}\chi_N(\langle a,\kappa^n_t\rangle)=\int_H\chi_N(\langle a,y\rangle)\ \mu(dy).
		\end{align*}
		But the right hand side is strictly positive for $N$ large enough. Hence $\langle a, \kappa_t^n\rangle$, $n\in\mathbb{N}$, is bounded in $\mathbb{R}$. But since for all $f\in C_b(\mathbb{R})$
		\begin{align*}
		\lim_{n\rightarrow\infty}f(\langle a,\kappa^n_t\rangle)=\int f(\langle a, y\rangle)\ \mu(dy),
		\end{align*}
		all accumulation points of $\langle a,\kappa^n_t\rangle$, $n\in\mathbb{N}$, must coincide. Consequently, 
		\begin{align*}
		\lim_{n\rightarrow\infty}\langle a, \kappa_t^n\rangle\ \text{exists for all }a\in H,
		\end{align*}
		as a linear functional in $a\in H$, so by the uniform boundedness principle must be continuous on $H$. Hence, there exists $\kappa_t\in H$ such that $\kappa_t^n\rightharpoonup\kappa_t$ in $H$ as $n\rightarrow\infty$. Taking $n\rightarrow\infty$ in $(\ref{R11})$, we therefore obtain that for all $a\in H$
		\begin{align*}
		\widehat{\nu_t}(a)=e^{i\langle a,\kappa_t\rangle}\ e^{-\int^t_{-\infty}\lambda(\sigma^*(r)U^*_{r,t}a)dr},
		\end{align*}
		i.e. 
		\begin{align*}
		\nu_t=\delta_{\kappa_t}\ast\mu_{-\infty,t}.
		\end{align*}
		So, for all $a\in H$
		\begin{align*}
		\int_H\lvert\langle a,y\rangle\rvert\ \nu_t(dy)\leq\lvert\langle a,\kappa_t\rangle\rvert+\int_H\lvert\langle a,y\rangle\rvert\ \mu_{-\infty,t}(dy)<\infty
		\end{align*}
		and 
		\begin{align*}
		\int_H\lvert\langle a,y\rangle\rvert\ \nu_t(dy)=\langle a,\kappa_t\rangle,
		\end{align*}
		therefore $\nu\in\mathcal{K}_e^1(\pi)$ and 
		\begin{align*}
		{\bf{p}}(\nu)=\kappa:=(\kappa_t)_{t\in\mathbb{R}}.
		\end{align*}
		By Proposition \ref{33!}, we have that $\kappa\in\mathcal{K}(U)$. Hence Claim {\emph{(ii)}} is proved.
		\begin{description}
			\item {\it{Claim (iii)}} {\emph{Let $\kappa\in\mathcal{K}(U)$. Then $\nu^{\kappa}\in\mathcal{K}^1_e(\pi)$.
					In particular, ${\bf{p}}\mid_{\mathcal{K}_e^1(\pi)}:\ \mathcal{K}_e^1(\pi)\rightarrow \mathcal{K}(U)$ is onto.}}
		\end{description}
		{\it{Proof}}: By Claim {\emph{(i)}} we have $\nu^\kappa\in\mathcal{K}^1(\pi)$, and thus by Lemma \ref{pp} and Claim {\emph{(ii)}}
		\begin{align*}
		\nu^\kappa&~=\int_{\mathcal{K}^1_e(\pi)}\widetilde{\nu}\ \xi_{\nu^\kappa}(d\widetilde{\nu})\\
		&~=\int_{\mathcal{K}^1_e(\pi)}\nu^{\bf{p}(\widetilde{\nu})}\ \xi_{\nu^\kappa}(d\widetilde{\nu})\\
		&~=\int_{\mathcal{K}(U)}\nu^{\widetilde{\kappa}}\ \eta(d\widetilde{\kappa}),
		\end{align*}
		where $\eta:=\xi_{\nu^\kappa}\circ{\bf{p}}^{-1}$, i.e. the image measure of $\xi_{\nu^\kappa}$ under ${\bf{p}}$ on $\mathcal{K}(U)\subset C(\mathbb{R};H)$ (see Remark \ref{R5} {\emph{(i)}}) equipped with the Borel $\sigma$-algebra inherited from $C(\mathbb{R};H)$ and where we have adapted the notation from the proof of Proposition \ref{R8}.\\ We claim that $\eta=\delta_{\kappa}$.\\
		Let $t\in\mathbb{R}$. Then for all $a\in H$
		\begin{align*}
		e^{i\langle a,\kappa_t\rangle}\cdot\widehat{\mu_{-\infty,t}}(a)
		=\widehat{\nu^\kappa_t}(a)=\int_{\mathcal{K}(U)}\widehat{\nu^{\widetilde{\kappa}}_t}(a)\ \eta(d\widetilde{\kappa})=\int_{\mathcal{K}(U)}e^{i\langle a,\widetilde{\kappa}_t\rangle}\widehat{\mu_{-\infty,t}}(a)\ \eta(d\widetilde{\kappa}).
		\end{align*}
		Since $\widehat{\mu_{-\infty,t}}(a)\neq0$ for any $a\in H$, we deduce that
		\begin{align*}
		\widehat{\delta_{\kappa_t}}(a)=e^{i\langle a,\kappa_t\rangle}=\int_{\mathcal{K}(U)}e^{i\langle a,\widetilde{\kappa}_t\rangle}\ \eta(d\widetilde{\kappa})
		=\int_{H}e^{i\langle a,h\rangle}\ (\eta\circ {\emph{pr}}_t^{-1})(dh)=\widehat{(\eta\circ {\emph{pr}}_t^{-1})}(a),
		\end{align*}
		where ${\emph{pr}}_t:\mathcal{K}(U)\rightarrow H$ with ${\emph{pr}}_t(\kappa)=\kappa_t$ for every $t\in\mathbb{R}$. Therefore, $\eta$ is a measure on $\mathcal{K}(U)$ such that
		$\delta_{\kappa_t}=\eta\circ {\emph{pr}}_t^{-1}$.\\
		For $t_1<\cdots<t_n$, let 
		${\emph{pr}_{t_1,\cdots,t_n}}:\mathcal{K}(U)\rightarrow H^{\{t_1,\cdots,t_n\}}$ denotes the map $(\kappa_t)_{t\in\mathbb{R}}\longmapsto(\kappa_{t_1},\cdots,\kappa_{t_n})$. As above it follows that  
		\begin{align*}
		\eta\circ{\emph{pr}^{-1}_{t_1,\cdots,t_n}}=\delta_{\kappa_{t_1}}\otimes\cdots\otimes\delta_{\kappa_{t_n}}.
		\end{align*}
		Then a monotone class argument implies that $\eta=\delta_{\kappa}$.\\
		Hence, also Claim {\emph{(iii)}} is proved.
	\end{proof}
\end{theorem}

\section{An application: uniqueness of the entrance law associated with $T$-periodic time-inhomogeneous (generalized) Mehler semigroups}\label{S4}
We recall that $U=(U_{s,t})_{s\leq t}$ is called $T$-periodic if $U_{s+T,t+T}=U_{s,t}$ for every $s\leq t$.
\begin{theorem}\label{3'}
	Assume that $(H1)-(H5)$ hold and that $U$ and $\sigma$ are $T$-periodic. Furthermore, suppose there exist $c,\omega\in(0,\infty)$ such that
	$ \lVert U(s,t)\rVert_{\mathcal{L}(H)}\leq c\ e^{-\omega(t-s)}$ for every $s\leq t$. Then, $(\mu_{-\infty,t})_{t\in\mathbb{R}}$ defined in $(H4)$ is the unique $T$-periodic $\pi$-entrance law in $\mathcal{K}^1(\pi)$.
	\begin{proof}
		Let $\nu\in\mathcal{K}^1(\pi)$, $\nu$ $T$-periodic. Then by Proposition \ref{R8} for all $a\in H$, $t\in\mathbb{R}$
		\begin{align}\label{R12}
		\widehat{\nu_t}(a)&~=\int_{\mathcal{K}_e^1(\pi)} e^{i\langle a,\widetilde{\kappa}_t\rangle}\ \eta(d\widetilde{\kappa})\ e^{-\int_{-\infty}^t\lambda(\sigma^*(r)U_{r,t}^*a)dr}\nonumber\\
		&~=\int_{H} e^{i\langle a,h\rangle}\eta_t(dh)\ \widehat{\mu_{-\infty,t}}(a),
		\end{align}
		where $\eta_t:=\eta\circ\ pr_t^{-1}$ and $\eta,\ pr_t$ are as defined in the proof of Claim {\emph{(iii)}} in the proof of Theorem \ref{27!}. Since $\nu_{t+T}=\nu_t$ and $\mu_{-\infty,t+T}=\mu_{-\infty,t}$ for all $t\in\mathbb{R}$, it follows from $(\ref{R12})$ that 
		\begin{align*}
		\widehat{\eta_{t+T}}(a)=\widehat{\nu_{t+T}}(a)\frac{1}{\widehat{\mu_{-\infty,t+T}}(a)}=\widehat{\nu_t}(a)\frac{1}{\widehat{\mu_{-\infty,t}}(a)}=\widehat{\eta_t}(a),\quad \forall a\in H.
		\end{align*}
		Hence $\eta_{t+T}=\eta_t$ for all $t\in\mathbb{R}$ and therefore by $(\ref{R12})$ for all $t\in\mathbb{R}$, $a\in H$ and $n\in\mathbb{N}$
		\begin{align}\label{R13}
		\widehat{\nu_t}(a)=\widehat{\eta_{t+nT}}(a)\ \widehat{\mu_{-\infty,t}}(a).
		\end{align}
		But by definition of $\eta_t$, we have for all $n\in\mathbb{N}$
		\begin{align*}
		\widehat{\eta_{t+nT}}(a)&~=\int_{\mathcal{K}_e^1(\pi)}e^{i\langle a,\widetilde{\kappa}_{t+nT}\rangle}\eta(d\widetilde{\kappa})\\
		&~=\int_{\mathcal{K}_e^1(\pi)}e^{i\langle a,U_{t,t+nT}\widetilde{\kappa}_t\rangle}\eta(d\widetilde{\kappa})\\
		&~=\int_{H}e^{i\langle a,U_{t-nT,t}h\rangle}\eta_t(dh)
		\end{align*}
		by the $T$-periodicity of $(U_{s,t})_{s\leq t}$.\\
		Hence by $(\ref{R13})$ and Lebesque's dominated convergence theorem for all $a\in H$, $t\in \mathbb{R}$
		\begin{align*}
		\widehat{\nu_t}(a)=\lim_{n\rightarrow\infty}\int_H e^{i\langle a,U_{t-nT,t}h\rangle}\eta_t(dh)\cdot\widehat{\mu_{-\infty,t}}(a)=\widehat{\mu_{-\infty,t}}(a),
		\end{align*}
		since $\lim_{n\rightarrow\infty}U_{t-nT,t}h=0$ for all $h\in H$. Therefore, $\nu_t=\mu_{-\infty,t}$ for all $t\in\mathbb{R}$ and Theorem \ref{3'} is proved.
	\end{proof} 
\end{theorem}
\begin{remark}
	For a related result under a different set of assumptions we refer to \cite[Theorem 4.11]{MR2861314}. Our proof is, however, considerably shorter than that in \cite{MR2861314}. 
	In the special Gaussian case (i.e. $M$ in $(\ref{101'})$ is the zero measure) the above theorem was first proved in \cite{MR2369672}.
\end{remark}

\section{Examples}
In this section, we are going to present two type of examples. First, we consider strong evolution families $(U_{s,t})_{s\leq t}$ as in $(H1)$ with bounded generators and a class of functions $\lambda$ as in $(H3)$, but being additionally Sazonov continuous. Second, we consider $(U_{s,t})_{s\leq t}$ with unbounded generators and a concrete $\lambda$ as in $(H3)$, which merely satisfies $(H4)$. In both cases, for simplicity we restrict to time homogeneous evolution families, but easy modifications then also lead to examples in the non-time homogeneous case.\\
So, let $H$ be a separable real Hilbert space as in the previous section and we fix $\sigma:\mathbb{R}\rightarrow\mathcal{L}(H)$ as in $(H2)$. We start with the following lemma, which will be very useful below. The proof is standard, but we include it for the reader's convenience.
\begin{lemma}\label{R19}
	Let $\vartheta$ be a finite positive measure on $(H,\mathscr{B}(H))$ and $\alpha\in[1,2]$ such that 
	\begin{align*}
	\int_H\lvert\langle a,x\rangle\rvert^\alpha\ \vartheta(dx)<\infty,\quad\forall a\in H.
	\end{align*}
	Then the map 
	\begin{align*}
	H\ni a\mapsto \int_H\lvert\langle a,x\rangle\rvert^\alpha\ \vartheta(dx) 
	\end{align*}
	is Sazonov continuous.
	\begin{proof}
		Since $\alpha\in[1,2]$, it obviously suffices to prove Sazonov continuity in $a=0$. So, let $\varepsilon\in(0,1)$ and $\mathcal{R}_\varepsilon\in(0,\infty)$ such that
		\begin{align*}
		\int_{\{\lVert\cdot\rVert>\mathcal{R}_\varepsilon\}}\lvert\langle a,x\rangle\rvert^\alpha\ \vartheta(dx)<\frac{\varepsilon}{2}.
		\end{align*}
		Recall that the covariance operator $\mathcal{S}_{\varepsilon}\in\mathcal{L}(H)$ defined by 
		\begin{align*}
		\int_{\{\lVert\cdot\rVert\leq\mathcal{R}_\varepsilon\}}\langle a_1,x\rangle\ \langle a_2,x\rangle\ \vartheta(dx)=\langle\mathcal{S}_\varepsilon a_1,a_2\rangle,\quad a_1,a_2\in H,
		\end{align*}	
		is symmetric, positive definite and of trace class. Hence, if 
		\begin{align*}
		a\in \big\{x\in H\mid\lVert\mathcal{S}^{\frac{1}{2}}_\varepsilon x\rVert<(\frac{\varepsilon}{2})^{\frac{2}{\alpha}}\vartheta(H)^{\frac{\alpha-2}{\alpha}}\big\},
		\end{align*}
		we have 
		\begin{align*}
		\int_H\lvert\langle a,x\rangle\rvert^\alpha\ \vartheta(dx)\leq\vartheta(H)^{\frac{2-\alpha}{2}}\bigg(\int_{\{\lVert\cdot\rVert\leq\mathcal{R}_\varepsilon\}}\langle a,x\rangle^2 \ \vartheta(dx)\bigg)^{\alpha/2}+\frac{\varepsilon}{2}<\varepsilon.
		\end{align*}
		Since $\mathcal{S}_\varepsilon^{\frac{1}{2}}$ is Hilbert-Schmidt, the assertion follows.
	\end{proof}
\end{lemma}
{\bf{\subsection{Bounded generators}}}
Let $\omega\in(0,\infty)$ and for $s,t\in\mathbb{R}$, $s\leq t$,
\begin{align}\label{R17}
U_{s,t}:=e^{-\omega(t-s)}I_H,
\end{align}
where $I_H$ denotes the identity map on $H$. Then obviously $(U_{s,t})_{s\leq t}$ is an evolution family satisfying $(H1)$ and $A(t):=-\omega\ e^{\omega t}I_H$, $t\in\mathbb{R}$, are the corresponding generators.
Furthermore, clearly $(U_{s,t})_{s\leq t}$ is strictly contractive, i.e. it satisfies condition $(U.1)$ in Remark \ref{R18} {\emph{(iii)}}.\\
We shall now define a class of $\lambda: H\rightarrow\mathbb{C}$ satisfying $(H3)$, $(H5)$ and $(\lambda.1)$, $(\lambda.2)$ in Remark \ref{R18} {\emph{(iii)}}, which hence by the latter satisfy $(H4)$ and our main result Theorem \ref{27!} applies to such $\lambda$ and $(U_{s,t})_{s\leq t}$ as in $(\ref{R17})$.\\
Let $\vartheta$ be as in Lemma \ref{R19} with $\alpha\in(1,2)$. Define
\begin{align}\label{R20}
\lambda(a):=\int_H\lvert\langle a,x\rangle\rvert^\alpha\ \vartheta(dx),\quad a\in H.
\end{align}
Since $s\mapsto\lvert s\rvert^\alpha$ is negative definite, $\lambda$ is negative definite. Therefore, since it is Sazonov continuous by Lemma \ref{R19}, hence norm-continuous, it clearly satisfies $(H3)$, $(H5)$ from Section 3 as well as $(\lambda.1)$ from Remark \ref{R18} ${\emph{(iii)}}$. So, it remains to prove $(\lambda.2)$.\\
To this end, we first note that by \cite{MR758255}, Proposition 6.4.5 and its proof, we know that 
\begin{align*}
\int_H\lVert x\rVert^\alpha\ d\vartheta(x)<\infty
\end{align*}
and that the L\'evy measure $M$ of $\lambda$ is given by 
\begin{align*}
M(B):= c_\alpha^{-1}\ \int_H\int_0^\infty\mathbbm{1}_B(tx)\ t^{-1-\alpha}\ dt\ \vartheta(dx),\quad B\in\mathscr{B}(H),
\end{align*}
where $c_\alpha\in(0,\infty)$. Hence
\begin{align*}
\int_{\{\lVert\cdot\rVert>1\}}\lVert x\rVert\ M(dx)&~=c_\alpha^{-1}\ \int_H\int_{\frac{1}{\lVert x\rVert}}^\infty\lVert x\rVert\ t^{-\alpha}dt\ \vartheta(dx)\\
&~=c_\alpha^{-1}\ (1-\alpha)^{-1}\ \int_H \lVert x\rVert^\alpha\ \vartheta(dx)<\infty.
\end{align*}
\begin{remark}\label{R23}
	{\emph{(i)}} We note that for $\alpha\in(1,2)$ and any symmetric, positive definite $\mathcal{S}\in\mathcal{L}(H)$ of trace class, the function
	\begin{align*}
	\lambda(a):=\lVert \mathcal{S}^{\frac{1}{2}}\ a\rVert^\alpha,\quad a\in H,
	\end{align*}
	is of type $(\ref{R20})$. Indeed, let $\mathcal{N}(0,\mathcal{S})$ be the centered Gaussian measure on $(H,\mathscr{B}(H))$ with covariance operator $\mathcal{S}$. Then an elementary calculation shows that for some constant $c_\alpha\in(0,\infty)$
	\begin{align*}
	\lambda(a)=c_\alpha\int_H\lvert\langle a,x\rangle \rvert^\alpha\ \mathcal{N}(0,\mathcal{S})(dx),\quad\forall a\in H.
	\end{align*}  
	{\emph{(ii)}} For our simple evolution family $(U_{s,t})_{s\leq t}$ defined in $(\ref{R17})$, we obviously have that 
	\begin{align*}
	\mathcal{K}(U)=\{\mathbb{R}\ni s\mapsto e^{-\omega s}x\mid x\in H\},
	\end{align*}
	i.e. $\mathcal{K}(U)$ is isomorphic to all of $H$.
\end{remark}
{\bf{\subsection{Unbounded generators}}}
Let $(A,\mathcal{D}(A))$ be a self-adjoint operator on $H$ such that for some $\omega\in(0,\infty)$
\begin{align*}
\langle Ax,x\rangle\leq -\omega\ \lVert x\rVert^2,\quad\forall x\in H,
\end{align*}
and that $(-A)^{-1}$ is trace class. Let $\{e_i\mid i\in\mathbb{N}\}$ be an eigenbasis of $A$ and $-\lambda_i$, $\lambda_i\in(0,\infty)$, be the corresponding eigenvalues, with $\lambda_i$ numbered in increasing order. Hence,
\begin{align}\label{R21}
\sum_{i=1}^\infty \frac{1}{\lambda_i}<\infty.
\end{align}  
\begin{example}
	Let $H:=L^2 \big((0,1),d\xi\big)$ with $d\xi=$Lebesgue measure and $A=\Delta$ with $\mathcal{D}(A):=H_0^1\big((0,1)\big)\cap H^2\big((0,1)\big)$, where the latter are the standard Sobolev spaces in $L^2\big((0,1),d\xi\big)$ of order 1 and 2, respectively, where the subscript zero refers to Dirichlet boundary conditions.
\end{example} 
Let
\begin{align*}
U_t:= e^{tA},\quad t>0.
\end{align*}
Then 
\begin{align}\label{R22}
U_{s,t}:= U_{t-s}=e^{(t-s)A},\quad s\leq t,
\end{align}
defines an evolution family satisfying $(H1)$ and $(U.1)$ from Remark \ref{R18} {\emph{(iii)}}.\\
Fix $\alpha\in(1,2)$ and define
\begin{align*}
\lambda(a):= \lVert a\rVert^\alpha,\quad a\in H.
\end{align*}
Then obviously $\lambda$ satisfies $(H3)$ and $(H5)$.\\
Now we are going to prove that $(H4)$ also holds:
Fix $t\in\mathbb{R}$. Then, we have by changing variables
\begin{align*}
\varPsi_t(a):=\int_{-\infty}^{t}\lambda\big(\sigma^*(r)\ U_{t-r}a\big)\ dr=\int_0^\infty \lambda\big(\sigma^*(t-r)\ U_{r}a\big)\ dr,\quad\forall a\in H,
\end{align*}
where the last integral is finite, since for $a\neq0$ by $(H2)$ it is up to a constant bounded by 
\begin{align*}
\int_0^\infty\bigg(\sum_{i=1}^\infty\ \langle U_ra,e_i\rangle^2\bigg)^{\alpha/2}dr&~=\int_0^\infty\bigg(\sum_{i=1}^\infty\ \langle a,e^{-\lambda_ir}e_i\rangle^2\bigg)^{\alpha/2}dr\\
&~=\lVert a \rVert^\alpha\ \int_0^\infty\ e^{-\omega\alpha r}dr\\
&~= \lVert a \rVert^\alpha\ \frac{1}{\omega\alpha}.
\end{align*}
We are now going to construct a finite measure $\vartheta$ on $(H,\mathscr{B}(H))$ such that 
\begin{align*}
\varPsi_t(a)=\int_H\lvert\langle a,x\rangle\rvert^\alpha\vartheta(dx),\quad \forall a\in H,
\end{align*}
which by Lemma \ref{R19} implies that $\varPsi_t$ is Sazonov continuous, which implies the first requirement in $(H4)$.\\
Clearly, by $(\ref{R21})$ also the linear operators $U_r=e^{rA}$, $r\in(0,\infty)$, are all symmetric, positive definite and of trace class, hence so are the operators 
\begin{align*}
\mathcal{S}_{r,t}:=\big(\rho(r)\big)^{-\frac{2}{\alpha}}\ U_r\ \sigma(t-s)\ \sigma^*(t-s)\ U_r,\quad r\in(0,\infty),
\end{align*} 
where $\rho\in L^1\big([0,\infty),dr\big)$ is a fixed function, $\rho>0$.\\
Therefore, for $r\in(0,\infty)$, we can consider $\mathcal{N}(0,\mathcal{S}_{r,t})$, i.e. the centered Gaussian measure on $(H,\mathscr{B}(H))$ with covariance operator $\mathcal{S}_{r,t}$. Then, as in Remark \ref{R23} {\emph{(i)}}
\begin{align*}
\lVert \sigma^*(t-r)\ U_ra\rVert^\alpha=c_\alpha\ \rho(r)\ \int_H\lvert\langle a,x\rangle\rvert^\alpha \ \mathcal{N}(0,\mathcal{S}_{r,t})(dx),\quad \forall a\in H
\end{align*}
for some $c_\alpha\in(0,\infty)$. Now define
\begin{align*}
\vartheta(dx):=c_\alpha\ \int_0^\infty\rho(r)\ \mathcal{N}(0,\mathcal{S}_{r,t})(dx)\ dr,
\end{align*}
which is a finite measure on $(H,\mathscr{B}(H))$ and we have
\begin{align*}
\varPsi_t(a)=\int_H\lvert\langle a,x \rangle\rvert^\alpha\ \vartheta(dx),\quad\forall a\in H.
\end{align*}
Hence, the measure $\mu_{-\infty,t}$ from $(H4)$ exists. It remains to show that it has weak first moments.\\
To this end, we again use \cite{MR758255}, Proposition 6.4.5 and its proof, to conclude that, for the L\'evy measure $M_t$ of $\varPsi_t$, we have
\begin{align}\label{R24}
\int_{\{\lVert\cdot\rVert>1\}}\ \lVert x\rVert\ M_t(dx)<\infty.
\end{align}
Let $a\in H$. Then since the L\'evy measure of $\mu_{-\infty,t}\circ\langle a,\cdot\rangle^{-1}$ is $M_t\circ\langle a,\cdot\rangle^{-1}$, by \cite{MR1739520}, Theorem 25.3, we only need to show that 
\begin{align}\label{R25}
\int_{\{\lvert\cdot\rvert>1\}}\ \lvert s\rvert\ \big(M_t\circ\langle a,\cdot\rangle^{-1}\big)(ds)<\infty.
\end{align} 
But the left hand side of $(\ref{R25})$ is equal to 
\begin{align*}
&~\int_{\{\lvert\langle a,\cdot\rangle\rvert>1\}}\ \lvert\langle a,x\rangle\rvert\ M_t(dx)\\\leq
&~\int_{\{\lVert \cdot\rVert\leq 1\}}\ \lvert\langle a,x\rangle\rvert^2\ M_t(dx)+\int_{\{\lVert \cdot\rVert\geq 1\}}\ \lvert\langle a,x\rangle\rvert\ M_t(dx)<\infty,
\end{align*}
since $M_t$ is a L\'evy measure and because of $(\ref{R24})$.

\begin{remark}
	For $U=(U_{s,t})_{s\leq t}$ defined in $(\ref{R22})$, the set $\mathcal{K}(U)$ seems difficult to describe explicitly. It is, however, again very big, because e.g. for every $i\in\mathbb{N}$
	\begin{align*}
	\kappa_s:=e^{-\lambda_i s}\ e_i,\quad s\in\mathbb{R},
	\end{align*}
	is obviously an element in $\mathcal{K}(U)$, and hence all linear combinations thereof.
\end{remark}

\thanks{Acknowledgment: The first author is grateful to Dr.~Tatiana Pasurek for numerous discussions. Financial support by the DFG
through SFB-701 is gratefully acknowledged.}


\begin{thebibliography}{10}

\bibitem{MR0481057}
C.~Berg and G.~Forst.
\newblock {\em Potential theory on locally compact abelian groups}.
\newblock Springer-Verlag, New York-Heidelberg, 1975.
\newblock Ergebnisse der Mathematik und ihrer Grenzgebiete, Band 87.

\bibitem{MR1392452}
V.~I. Bogachev, M.~R\"{o}ckner, and B.~Schmuland.
\newblock Generalized {M}ehler semigroups and applications.
\newblock {\em Probab. Theory Related Fields}, 105(2):193--225, 1996.

\bibitem{MR2369672}
G.~Da~Prato and A.~Lunardi.
\newblock Ornstein-{U}hlenbeck operators with time periodic coefficients.
\newblock {\em J. Evol. Equ.}, 7(4):587--614, 2007.

\bibitem{MR2243880}
D.~A. Dawson and Z.~Li.
\newblock Skew convolution semigroups and affine {M}arkov processes.
\newblock {\em Ann. Probab.}, 34(3):1103--1142, 2006.

\bibitem{MR2048508}
D.~A. Dawson, Z.~Li, B.~Schmuland, and W.~Sun.
\newblock Generalized {M}ehler semigroups and catalytic branching processes
  with immigration.
\newblock {\em Potential Anal.}, 21(1):75--97, 2004.

\bibitem{MR0518321}
E.~B. Dynkin.
\newblock Sufficient statistics and extreme points.
\newblock {\em Ann. Probab.}, 6(5):705--730, 1978.

\bibitem{MR1013934}
E.~B. Dynkin.
\newblock Three classes of infinite-dimensional diffusions.
\newblock {\em J. Funct. Anal.}, 86(1):75--110, 1989.

\bibitem{MR1745332}
M.~Fuhrman and M.~R\"{o}ckner.
\newblock Generalized {M}ehler semigroups: the non-{G}aussian case.
\newblock {\em Potential Anal.}, 12(1):1--47, 2000.

\bibitem{MR2861314}
F.~Kn\"{a}ble.
\newblock Ornstein-{U}hlenbeck equations with time-dependent coefficients and
  {L}\'{e}vy noise in finite and infinite dimensions.
\newblock {\em J. Evol. Equ.}, 11(4):959--993, 2011.

\bibitem{MR1930955}
P.~Lescot and M.~R\"{o}ckner.
\newblock Generators of {M}ehler-type semigroups as pseudo-differential
  operators.
\newblock {\em Infin. Dimens. Anal. Quantum Probab. Relat. Top.},
  5(3):297--315, 2002.

\bibitem{MR2032114}
P.~Lescot and M.~R\"{o}ckner.
\newblock Perturbations of generalized {M}ehler semigroups and applications to
  stochastic heat equations with {L}evy noise and singular drift.
\newblock {\em Potential Anal.}, 20(4):317--344, 2004.

\bibitem{MR2760602}
Z.~Li.
\newblock {\em Measure-valued branching {M}arkov processes}.
\newblock Probability and its Applications (New York). Springer, Heidelberg,
  2011.

\bibitem{MR2105913}
Z.~Li and Z.~Wang.
\newblock Generalized {M}ehler semigroups and {O}rnstein-{U}hlenbeck processes
  arising from superprocesses over the real line.
\newblock {\em Infin. Dimens. Anal. Quantum Probab. Relat. Top.},
  7(4):591--605, 2004.

\bibitem{MR758255}
W.~Linde.
\newblock {\em Infinitely divisible and stable measures on {B}anach spaces},
  volume~58 of {\em Teubner-Texte zur Mathematik [Teubner Texts in
  Mathematics]}.
\newblock BSB B. G. Teubner Verlagsgesellschaft, Leipzig, 1983.
\newblock With German, French and Russian summaries.

\bibitem{Lunardi2019SchauderTF}
A.~Lunardi and M.~R\"{o}ckner.
\newblock Schauder theorems for a class of (pseudo-)differential operators on
  finite and infinite dimensional state spaces, 2019.

\bibitem{MR3466574}
S.-X. Ouyang and M.~R\"{o}ckner.
\newblock Time inhomogeneous generalized {M}ehler semigroups and skew
  convolution equations.
\newblock {\em Forum Math.}, 28(2):339--376, 2016.

\bibitem{MR2886463}
S.-X. Ouyang, M.~R\"{o}ckner, and F.-Y. Wang.
\newblock Harnack inequalities and applications for {O}rnstein--{U}hlenbeck
  semigroups with jump.
\newblock {\em Potential Anal.}, 36(2):301--315, 2012.

\bibitem{MR0226684}
K.~R. Parthasarathy.
\newblock {\em Probability measures on metric spaces}.
\newblock Probability and Mathematical Statistics, No. 3. Academic Press, Inc.,
  New York-London, 1967.

\bibitem{MR1159586}
M.~R\"{o}ckner.
\newblock On the parabolic {M}artin boundary of the {O}rnstein-{U}hlenbeck
  operator on {W}iener space.
\newblock {\em Ann. Probab.}, 20(2):1063--1085, 1992.

\bibitem{MR1996872}
M.~R\"{o}ckner and F.-Y. Wang.
\newblock Harnack and functional inequalities for generalized {M}ehler
  semigroups.
\newblock {\em J. Funct. Anal.}, 203(1):237--261, 2003.

\bibitem{MR1739520}
K.-i. Sato.
\newblock {\em L\'{e}vy processes and infinitely divisible distributions},
  volume~68 of {\em Cambridge Studies in Advanced Mathematics}.
\newblock Cambridge University Press, Cambridge, 1999.
\newblock Translated from the 1990 Japanese original, Revised by the author.

\bibitem{MR1841407}
B.~Schmuland and W.~Sun.
\newblock On the equation {$\mu_{t+s}=\mu_s*T_s\mu_t$}.
\newblock {\em Statist. Probab. Lett.}, 52(2):183--188, 2001.

\bibitem{MR1435288}
N.~N. Vakhania, V.~I. Tarieladze, and S.~A. Chobanyan.
\newblock {\em Probability distributions on {B}anach spaces}, volume~14 of {\em
  Mathematics and its Applications (Soviet Series)}.
\newblock D. Reidel Publishing Co., Dordrecht, 1987.
\newblock Translated from the Russian and with a preface by Wojbor A.
  Woyczynski.

\end{thebibliography}
\end{document}